\documentclass[12pt]{amsart}
\usepackage{graphicx, amscd}
\usepackage[mathscr]{eucal}

\theoremstyle{plain}
\newtheorem{theorem}{Theorem}

\newtheorem{proposition}{Proposition}
\newtheorem{lemma}{Lemma}

\theoremstyle{definition}
\newtheorem{remark}{Remark}
\newtheorem{demo}{Proof}

\def\UU{\mathcal {U}}
\def\CC{\mathcal {C}}
\def\XX{\mathcal {X}}
\def\ZZ{\mathbb {Z}}
\def\C{\mathbb {C}}
\def\PP{\mathbb {P}}
\def\AA{\mathbb {A}}
\DeclareMathOperator\Spec{Spec}

\begin{document}
\author{Sergey Galkin}
\address{HSE University, Russian Federation}
\author{Sergey Rybakov}
\footnotetext{This work is supported by the Russian Science Foundation under grant № 18-11-00121}
\address{Institute for information transmission problems of the Russian Academy of Sciences, AG Laboratory, HSE University, Russian Federation}

\email{rybakov.sergey@gmail.com}%
\title{A family of K3 surfaces and towers of algebraic curves over finite fields}
\date{}


\maketitle


\section{Introduction}
Let $\bar k$ be an algebraic closure of a finite field $k$ of characteristic $p$,
and $\overline{Y}$ be the base change $$\overline{Y} = Y\times_{\Spec k}{\Spec\bar k}$$ of an algebraic 
variety $Y$ over $k$. We also fix a prime number $\ell\neq p$. 

In the previous paper~\cite{Ry} we constructed an optimal tower of algebraic curves over ${\mathbb F}_{p^2}$ starting from a 
smooth family $f:X\to U$ of elliptic curves over an open subset $U\subset{\mathbb P}^1$.
We used that the $i$-th derived \'etale direct image of the constant sheaf ${\mathbb Z}/\ell^n {\mathbb Z}$ 
corresponds to a local system ${\mathcal{V}}_n$ on $U$, and that 
one can define a fiberwise \emph{projectivisation} $P_n({\mathcal{V}}_n)$
of this local system, which is an \'etale scheme $U_n$ over $U$ (see section~\ref{s2}).  
The challenge of this approach is to bound the genus of a connected component of $U_n$.

In this paper we use a special case of this construction, that lifts to characteristic zero; thus we
 can bound the ramification and the genus in the tower using transcendental methods.
We give an example of a good tower coming from a family $V_\lambda$ of $K3$ surfaces (see section~\ref{s3}).
The family that we consider is the mirror-dual family to $\PP^3$.
It was studied, among others, by Dolgachev~\cite{Dolg}, Batyrev~\cite{Bat},
Coates--Corti--Galkin--Kasprzyk~\cite{CCGK}. 
Namely, Batyrev proved that $V_\lambda$ is mirror-dual to $\PP^3$,
and Dolgachev computed global monodromy of $V_\lambda$. 
We also use explicit description of local monodromy operators given in~\cite{PS}.
Altogether this information is enough to describe the corresponding tower.

\begin{theorem}\label{T1}
There exists a tower $C_\bullet$ over $k=\mathbb{F}_{p^2}$ that is optimal for $p=3$, and is good if $p\equiv 3(4)$.
\end{theorem}

Other similar families that can be used as an input data for our construction
can be found in~\cite{Gol} and~\cite{CCGK}. 

\section{Towers of algebraic curves}\label{s2}
In this section we give a modified version of the construction of a tower of algebraic curves from~\cite{Ry}.

First we recall a generalization of the projectivization of an ${\mathbb F}_\ell$-vector space.
Put $\Lambda_n={\mathbb Z}/\ell^n{\mathbb Z}$.
Let $V$ be a finitely generated $\Lambda_n$--module. The set 
$$V^*=\{v\in V|\ell^{n-1}v\neq 0\}$$ has a natural action of the group of 
invertible elements $\Lambda^*_n$. We say that the set $$P_n(V)=V^*/\Lambda_n^*$$ 
is \emph{the projectivisation of $V$}.  

\begin{lemma}\cite[Lemma 1]{Ry}\label{lem1}
Suppose that $V_n$ is a $\Lambda_n$--module, and $V_{n-1}$ is a $\Lambda_{n-1}$--module. 
The group $V_{n-1}$ is a $\Lambda_n$--module in an 
obvious way. Let $\varphi:V_n\to V_{n-1}$ be a homomorphism such that $\ell\ker\varphi=0$.
Then there is an induced map $$P(\varphi): P_n(V_n)\to P_{n-1}(V_{n-1}).$$ 
\end{lemma}

Let $\UU$ be a smooth and connected scheme over $\ZZ$ of relative dimension one, and let $f:\XX\to \UU$ be a smooth family of algebraic varieties over $\UU$. 
Suppose that there is an action of a finite group $H$ on $\XX$, and that $f$ is $H$-equivariant.
According to~\cite[VI.4.2]{Milne}, the $H$-invariants of the \'etale direct image
$${\mathcal{V}}_n=(R_{{\acute e}t}^if_{U*}(\Lambda_n))^H,$$
is a locally constant finite sheaf of $\Lambda_n$--modules on $\UU$. 

According to~\cite[V.1.1]{Milne}, the functor ${\mathcal{V}}_n^*$ given by 
$$T\mapsto {\mathcal{V}}_n(T)^*$$ is representable by a finite \'etale scheme over $\UU$ 
with a free action of the group $\Lambda_n^*$. 
Thus the quotient  $\widetilde{\UU}_n$ by this action represents the sheaf $P_n({\mathcal{V}}_n)$. 
For each $n$ we have a natural \'etale morphism $$\pi_n:\widetilde{\UU}_n\to \widetilde{\UU}_{n-1}.$$

Let $K$ be a field, and let $\Bar a\in \UU({\bar K})$ be a point of $\UU$.
The stalk $({\mathcal{V}}_n)_{\Bar a}$ of ${\mathcal{V}}_n$ is isomorphic to $H^i_{{\acute e}t}({\overline{X}}_{\Bar a},\Lambda_n)^H$. 
According to~\cite[Lemma 2]{Ry} the fiber of $\widetilde{\UU}_n$ over ${\bar a}\in \UU({\bar K})$ is naturally isomorphic to $P_n(({\mathcal{V}}_n)_{\bar a})$.
In what follows we assume that the module $H^i_{{\acute e}t}({\overline{X}}_{\Bar a},{\mathbb Z}_\ell)^H$ is free.

Let $\Bar a\in \UU(\C)$, and let $\pi_1(\UU,\Bar a)$ be the fundamental group of $\UU$. 

\begin{proposition}\label{prop2} 
There is a one to one correspondence between orbits of the action of $\pi_1(\UU,\Bar a)$ on  
$P_n(H^i_{{\acute e}t}(\overline{X}_{\bar a},\Lambda_n))$ and connected components of $\widetilde{\UU}_n$.
The degree of the connected component over $\UU$ is equal the cardinality of the corresponding orbit.
\end{proposition}
\begin{demo}
Let $K$ be a field of functions of $\UU$.
A connected component ${\mathcal {L}}$ of $\widetilde{\UU}_n$ is uniquely determined by its field of functions.
The group $\pi_1(\UU,{\bar a})$ is isomorphic to the quotient of the Galois group ${\mathrm {Gal}}(K)$. 
Finite field extensions of $K$ correspond to transitive actions of ${\mathrm {Gal}}(K)$ on finite sets (see section 2.3 of ~\cite{Ry}). 
Finally, the degree of ${\mathcal {L}}$ over $\UU$ is equal to the degree of the corresponding field extension.
\end{demo}

We assume that $U=\UU\otimes k$ is non-empty and geometrically connected. 
If $\UU_n$ is a connected component of $\widetilde{\UU}_n$, then $\UU_{n-1}=\pi_n(\UU_n)$ is a connected component of $\widetilde{\UU}_n$.
Choose a sequence of connected subschemes $\UU_n\subset\widetilde{\UU}_n$ such that $\UU_{n-1}=\pi_n(\UU_n)$.
Let $C=C_0$ be a smooth projectivisation of $\UU\otimes k$, and let $C_n$ be a regular projectivisation of $\UU_n\otimes k$. 
Note that every connected component $Y$ of $C_n$ is a smooth algebraic curve over the algebraic closure of $k$ in $k(Y)$.
Put $W_n=\UU_n\otimes\C$. The smooth projectivisation $\CC_n$ of $W_n$ is an algebraic curve over $\C$. We get the following diagram:

$$\begin{array}{ccccccc}
\downarrow &  &\downarrow & \downarrow &  \downarrow &  &\downarrow\\
C_n &\supset & U_n & \UU_n & W_n & \subset & \CC_n\\
\downarrow &  &\downarrow & \downarrow &  \downarrow &  &\downarrow\\
\vdots &  &\vdots & \vdots &  \vdots &  &\vdots\\
\downarrow &  &\downarrow & \downarrow &  \downarrow &  &\downarrow\\
C_1 &\supset & U_1 & \UU_1 & W_1 & \subset & \CC_1\\
\downarrow &  &\downarrow & \downarrow &  \downarrow &  &\downarrow\\
C &\supset & U & \UU & W & \subset & \CC\\
\end{array}
$$


Let $a\in U(k)$. 
We say that a smooth variety $Y$ over $k$ is \emph{strongly supersingular in degree $i$} if the Frobenius action on $H^i_{{\acute e}t}(\overline{Y},{\mathbb Z}_\ell)$ is the multiplication by 
$q^{i/2}$ or $-q^{i/2}$. 

\begin{proposition}\label{prop3}
Assume that $C_n$ is connected, and there exists $a\in U(k)$ such that $X_a$ is strongly supersingular in degree $i$. 
Then $C_n$ is a smooth, geometrically irreducible curve over ${\mathbb F}_p$, and the fiber of the morphism $C_n\to C_0$ over $a$ is split over $k$.
\end{proposition}
\begin{demo}
According to~\cite[Corollary 3]{Ry} the statement is true after the base change to $k$; thus $C_n$ is smooth and geometrically irreducible over ${\mathbb F}_p$ as well. 
\end{demo}

The ramification of $\CC_n$ over $\CC$ can be computed in terms of the local monodromy of the family $f\otimes\C$. 
Fix a point $\Bar a\in W(\C)$. Let $y\in \CC(\C)-W(\C)$, and let $G_y$ be the stabilizator of $y$ in $\pi_1(W,\Bar a)$.

\begin{proposition}\label{prop4}\cite[Proposition 1]{Ry}
There is a bijection between orbits of the action of $G_y$ on
 $P_n(({\mathcal{V}}_n)_{\Bar a})$ and points $x\in \CC_n(\C)$ over $y$. 
The ramification index of $x$ is equal to the cardinality of the corresponding orbit.
\end{proposition}

\begin{remark}
Let $k$ be an extension of odd degree of ${\mathbb F}_p$.
Assume that the general fiber $Y$ of our family is a K3 surface. Then, by a result of Artin~\cite[(6.8)]{Ar}, there are no strongly supersingular fibers in our family.
This observation force us to look for a weaker supersingularity condition such that one still could prove an analog of Proposition~\ref{prop3} for $k=\mathbb{F}_p$.
\end{remark}

\section{A family of $K3$ surfaces}\label{s3}
In this section we prove Theorem~\ref{T1}.

Projection $\pi$ from
\[ X = \{x_0^4+x_1^4+x_2^4+x_3^4+4\lambda x_0x_1x_2x_3=0\}\subset \AA^1_\lambda\times\PP^3_x \]
to $\AA^1_\lambda$ is a pencil of quartic surfaces $X_\lambda$,
and its restriction to $\UU=\Spec\ZZ[\lambda,\frac{1}{2(\lambda^4-1)}]$ is smooth.

Let $i$ be a primitive $4$th root of unity, considered as a generator of a group scheme $\mu_4 = \Spec\ZZ[i]/(i^4-1)$.
Maps
\[(x_0:x_1:x_2:x_3)\to(x_0:ix_1:x_2:-ix_3)
\text{ and }
(x_0:x_1:x_2:x_3)\to(ix_0:x_1:-ix_2:x_3)\]
generate a non-effective action of a group scheme $H=\mu_4\times\mu_4$ on $X$.
The action of $H$ commutes with the projection $\pi$,
hence it induces an $H$-action on every fiber.

\begin{theorem}\cite[Theorem 7.1, Theorem 7.7, Theorem 8.2]{Dolg}
The general fiber of the quotient family $(X_\lambda\otimes \C)/H$ has $6$ rational double points of type $A_3$. 
The resolution of singularities is the family $V_\lambda$ of Kummer surfaces associated to the product $E_1\times E_2$, where
$E_1$ and $E_2$ are elliptic curves with $2$-isogeny between them. The monodromy group $G$ of this family is the modular Fricke group of level $2$, i.e.,
$G$ is the subgroup of $\mathrm{PGL}(2,\mathbb{R})$ generated by $\Gamma_0(2)$ and by
$$\left (\begin{matrix}
0 & -1/\sqrt{2}\\
\sqrt{2} & 0\\	
\end{matrix}\right).$$
\end{theorem}

Exceptional lines of the resolution of $6$ rational double points of type $A_3$ generate a sublattice of rank $18$ in $\mathrm{NS}(V_\lambda)$.
It follows that the rank of $H$-invariants $H^2(X_\lambda\otimes \C,\ZZ)^H$ is $4$. Moreover, a hyperplane section gives a $G$-invariant sublattice of rank $1$. 
The orthogonal complement $T$
is naturally isomorphic to the 
lattice $T$ of quadratic forms $$\{\alpha x^2+2\sqrt{2}\beta xy+\gamma y^2| \alpha,\beta,\gamma\in\ZZ\}$$ endowed with the quadratic form
 $Q(\alpha,\beta,\gamma)=2(2\beta^2-\alpha\gamma)$~\cite[Theorem 7.1]{Dolg}.
The action of $G$ on $T$ is induced from the following twisted action $\rho$ of $\mathrm{PGL}(2,\mathbb{R})$ on $T\otimes\mathbb{R}$: 
$$\rho\left (\begin{matrix}
 a& b\\
c & d\\	
\end{matrix}\right) =
\left (\begin{matrix}
 1/\sqrt[4]{2}& 0\\
0 & \sqrt[4]{2}\\	
\end{matrix}\right)
\left (\begin{matrix}
 a& b\\
c & d\\	
\end{matrix}\right)
\left (\begin{matrix}
 1/\sqrt[4]{2}& 0\\
0 & \sqrt[4]{2}\\	
\end{matrix}\right)^{-1}.$$

The zero locus of the quadratic form on $T$ is exactly the set of degenerate forms. Let $D_n\subset P_n(T/\ell^nT)$ be the projectivisation of the cone of degenerate forms.
It is straightforward to check that the action of $G$ on $D_n$ is transitive for all $n$. 
It follows from Proposition~\ref{prop2} that the curves $\CC_n$ are connected; therefore $C_n$ are connected as well.
Clearly, $D_1$ is a set of points on a conic in $P_1(T/\ell T)$; thus the degree of $\CC_1$ over $\CC_0$ is equal to $\ell+1$.
If $n>1$, then $\pi_n$ induces a morphism $\CC_n\to \CC_{n-1}$ of degree $\ell$.

\begin{lemma}
For the genus $g_n$ of $\CC_n$ we have 
$$g_n\leq (\ell+1)\ell^{n-1}/2.$$
\end{lemma}
\begin{demo}
The family $X_\lambda$ is mirror-dual to a family $V_\lambda$~\cite[Corollary 5.5.6]{Bat}.
In turn, the singular fibers and monodromy operators for such families are known.
According to~\cite[Period sequence 12]{PS} the mirror-dual family to $\PP^3$ has $4$ singular fibers such that the monodromy operator is an involution,
 and one fiber with maximal unipotent monodromy. The Hurwitz formula now gives
$$2g_n-2=-2(\ell+1)\ell^{n-1}+R_0+R,$$
where $R_0\leq \deg\pi_n=(\ell+1)\ell^{n-1}$ is the ramification degree of the point with maximal unipotent monodromy,
and $R$ is the ramification degree in other points.
We now bound $R$ using Proposition~\ref{prop4}. Since the monodromy operator is an involution, there are at least $(\ell+1)\ell^{n-1}/2$ orbits of degree $2$ on the set $D_n$.
This gives $$R\leq 4((\ell+1)\ell^{n-1}/2)=2(\ell+1)\ell^{n-1}.$$ This proves the lemma.
\end{demo}

By a result of Tsfasman~\cite[Example 5.9]{Ts}, the quartic $$x_0^4+x_1^4+x_2^4+x_3^4=0$$ is strongly supersingular in degree $2$ over $k$, when  $p\equiv 3(4)$. 
According to Proposition~\ref{prop3}, there exists at least one split fiber of $\pi_n$ over $k$. 
We finally get $$\lim_{n\to\infty}\frac{|C_n(k)|}{g_n}\geq 2.$$ This proves Theorem~\ref{T1}.

\end{document}